\documentclass[12pt,a4paper]{article}
\usepackage[utf8]{inputenc}
\usepackage[T1]{fontenc}
\usepackage[english]{babel}
\usepackage{amsmath, amssymb, amsthm}
\usepackage{geometry}
\usepackage{hyperref}
\newcommand{\eql}[1]{\begin{equation}#1\end{equation}}
\newcommand{\guill}[1]{``#1''}

\title{The Mathematics of Artificial Intelligence}
\author{
    Gabriel Peyr\'e \\
    CNRS and ENS, Universit\'e PSL \\
    \texttt{gabriel.peyre@ens.fr}
}
\date{\today}

\begin{document}

\maketitle

\begin{abstract}
This overview article highlights the critical role of mathematics in artificial intelligence (AI), emphasizing that mathematics provides tools to better understand and enhance AI systems. Conversely, AI raises new problems and drives the development of new mathematics at the intersection of various fields. This article focuses on the application of analytical and probabilistic tools to model neural network architectures and better understand their optimization. Statistical questions (particularly the generalization capacity of these networks) are intentionally set aside, though they are of crucial importance. We also shed light on the evolution of ideas that have enabled significant advances in AI through architectures tailored to specific tasks, each echoing distinct mathematical techniques. The goal is to encourage more mathematicians to take an interest in and contribute to this exciting field.
\end{abstract}

\section{Supervised Learning}

Recent advancements in artificial intelligence have mainly stemmed from the development of neural networks, particularly deep networks. The first significant successes, after 2010, came from supervised learning, where training pairs $(x^i, y^i)_i$ are provided, with $x^i$ representing data (e.g., images or text) and $y^i$ corresponding labels (typically, classes describing the content of the data). More recently, spectacular progress has been achieved in unsupervised learning, where labels $y_i$ are unavailable, thanks to techniques known as generative AI. These methods will be discussed in Section~\ref{sec:ia-gen}.

\paragraph{Empirical Risk Minimization.}
In supervised learning, the objective is to construct a function $f_\theta(x)$, dependent on parameters $\theta$, such that it approximates the data well: $y^i \approx f_\theta(x^i)$. The dominant paradigm involves finding the parameters $\theta$ by minimizing an empirical risk function, defined as:
\eql{\label{eq:erm}
\min_\theta E(\theta) := \sum_i \ell(f_\theta(x^i), y^i),
}
where $\ell$ is a loss function, typically $\ell(y, y') = \|y - y'\|^2$ for vector-valued $(y^i)_i$. Optimization of $E(\theta)$ is performed using gradient descent:
\eql{
\theta_{t+1} = \theta_t - \tau \nabla E(\theta_t), \label{eq:grad-desc}
}
where $\tau$ is the step size. In practice, a variant known as stochastic gradient descent \cite{robbins1951stochastic} is used to handle large datasets by randomly sampling a subset at each step $t$. A comprehensive theory for the convergence of this method exists in cases where $f_\theta$ depends linearly on $\theta$, as $E(\theta)$ is convex. However, in the case of deep neural networks, $E(\theta)$ is non-convex, making its theoretical analysis challenging and still largely unresolved. In some instances, partial mathematical analyses provide insights into the observed practical successes and guide necessary modifications to improve convergence. For a detailed overview of existing results, refer to \cite{bach2024learning}.

\paragraph{Automatic Differentiation.}
Computing the gradient $\nabla E(\theta)$ is fundamental. Finite difference methods or traditional differentiation of composite functions would be too costly, with complexity on the order of $O(D T)$, where $D$ is the dimensionality of $\theta$ and $T$ is the computation time of $f_\theta$.
The advances in deep network optimization rely on the use of automatic differentiation in reverse mode \cite{griewank2008evaluating}, often referred to as ``backpropagation of the gradient,'' with complexity on the order of $O(T)$.

\section{Two-Layer Neural Networks}

\paragraph{Multi-Layer Perceptrons.}
A multi-layer network, or multi-layer perceptron (MLP)~\cite{rosenblatt1958perceptron}, computes the function $f_\theta(x) = x_L$ in $L$ steps (or layers) recursively starting from $x_0=x$:
\eql{
x_{\ell+1} = \sigma(W_\ell x_\ell + b_\ell), \label{eq:mlp-multicouches}
}
where $W_\ell$ is a weight matrix, $b_\ell$ a bias vector, and $\sigma$ a non-linear activation function, such as the sigmoid $\sigma(s) = \frac{e^s}{1 + e^s}$, which is bounded, or ReLU (Rectified Linear Unit), $\sigma(s) = \max(s, 0)$, which is unbounded. The parameters to optimize are the weights and biases $\theta = (W_\ell,b_\ell)_{\ell}$. If $\sigma$ were linear, then $f_\theta$ would remain linear regardless of the number of layers. The non-linearity of $\sigma$ enriches the set of functions representable by $f_\theta$, by increasing both the width (dimension of intermediate vectors $x_k$) and the depth (number $L$ of layers).

\paragraph{Two-Layer Perceptrons and Universality.}
The mathematically best-understood case is that of $L = 2$ layers. Denoting $W_1 = (v_k)_{k=1}^n$ and $W_2^\top = (u_k)_{k=1}^n$ as the $n$ rows and columns of the two matrices (where $n$ is the network width), we can write $f_\theta$ as a sum of contributions from its $n$ neurons:
\eql{
    f_\theta(x) = \frac{1}{n} \sum_{k=1}^n u_k \sigma( \langle v_k, x \rangle + b_k ). \label{eq:mlp-2-couches}
}
The parameters are $\theta = (u_k, v_k, b_k)_{k=1}^n$. Here, we added a normalization factor $1/n$ (to later study the case $n \to +\infty$) and ignored the non-linearity of the second layer. A classical result by Cybenko~\cite{cybenko1989approximation} shows that these functions $f_\theta$ can uniformly approximate any continuous function on a compact domain. This result is similar to the Weierstrass theorem for polynomial approximation, except that \eqref{eq:mlp-2-couches} defines, for a fixed $n$, a non-linear function space. The elegant proof by Hornik~\cite{hornik1989multilayer} relies on the Stone-Weierstrass theorem, which implies the result for $\sigma = \cos$ (since \eqref{eq:mlp-2-couches} defines an algebra of functions when $n$ is arbitrary). The proof is then completed by uniformly approximating a cosine in 1-D using a sum of sigmoids or piecewise affine functions (e.g., for $\sigma = \text{ReLU}$).

\paragraph{Mean-Field Representation.}
This result is, however, disappointing, as it does not specify how many neurons are needed to achieve a given approximation error. This is impossible without adding assumptions about the regularity of the function to approximate. Barron's fundamental result~\cite{barron1993universal} introduces a Banach space of regular functions, defined by the semi-norm $\|f\|_B := \int |\hat{f}(\omega)| |\omega| \, \mathrm{d}\omega$, where $\hat{f}$ denotes the Fourier transform of $f$. Barron shows that, in the $L^2$ norm, for any $n$, there exists a network $f_\theta$ with $n$ neurons such that the approximation error on a compact $\Omega$ of radius $R$ is of the order $\|f - f_\theta\|_{L^2(\Omega)} = O\left(\frac{R \|f\|_B}{\sqrt{n}}\right)$. This result is remarkable because it avoids the \guill{curse of dimensionality}: unlike polynomial approximation, the error does not grow exponentially with the dimension $d$ (although the dimension affects the constant $\|f\|_B$).

Barron's proof relies on a \guill{mean-field} generalization of \eqref{eq:mlp-2-couches}, where a distribution (a probability measure) $\rho$ is considered over the parameter space $(u, v, b)$, and the network is expressed as:
\eql{
    F_\rho(x) := \int u \sigma(\langle v, x \rangle + b) \, \mathrm{d}\rho(u, v, b). \label{eq:mean-field}
}
A finite-size network \eqref{eq:mlp-2-couches}, $f_\theta = F_{\hat{\rho}}$, is obtained with a discrete measure $\hat{\rho} = \frac{1}{n} \sum_{k=1}^n \delta_{(u_k, v_k, b_k)}$.
An advantage of the representation~\eqref{eq:mean-field} is that $F_\rho$ depends linearly on $\rho$, which is crucial for Barron's proof. This proof uses a probabilistic Monte Carlo-like method (where the error decreases as $1/\sqrt{n}$ as desired): it involves constructing a distribution $\rho$ from $f$, then sampling a discrete measure $\hat{\rho}$ whose parameters $(u_k, v_k, b_k)_k$ are distributed according to $\rho$.

\paragraph{Wasserstein Gradient Flow.}
In general, analyzing the convergence of the optimization \eqref{eq:grad-desc} is challenging because the function $E$ is non-convex. Recent analyses have shown that when the number of neurons $n$ is large, the dynamics are not trapped in a local minimum. The fundamental result by Chizat and Bach~\cite{chizat2018global} is based on the fact that the distribution $\rho_t$ defined by gradient descent \eqref{eq:grad-desc}, as $\tau_t \to 0$, follows a gradient flow in the space of probability distributions equipped with the optimal transport distance. These Wasserstein gradient flows, introduced by~\cite{jordan1998variational} and studied extensively in the book~\cite{ambrosio2008gradient}, satisfy a McKean-Vlasov-type partial differential equation:
\eql{
    \partial_t \rho_t + \mathrm{div}(\rho_t \mathcal{V}(\rho)) = 0, \label{eq:edp-mlps}
}
where $x \to \mathcal{V}(\rho)(x)$ is a vector field depending on $\rho$, which can be computed explicitly from the data $(x^i, y^i)_i$. The result by Chizat and Bach can be viewed both as a PDE result (convergence to an equilibrium of a class of PDEs with a specific vector field $\mathcal{V}(\rho)$) and as a machine learning result (successful training of a two-layer network via gradient descent when the number of neurons is sufficiently large).

\section{Very Deep Networks}

The unprecedented recent success of neural networks began with the work of \cite{krizhevsky2012imagenet}, which demonstrated that deep networks, when trained on large datasets, achieve unmatched performance. The first key point is that to achieve these performances, it is necessary to use weight matrices $W_k$ that exploit the structure of the data. For images, this means using convolutions~\cite{lecun1998gradient}. However, this approach is not sufficient for extremely deep networks, with the number of layers $L$ reaching into the hundreds. 

\paragraph{Residual Networks.}
The major breakthrough that empirically demonstrated that network performance increases with $L$ was the introduction of residual connections, known as ResNet \cite{he2016deep}. The main idea is to ensure that, for most layers $x_\ell$, the dimensions of $x_\ell$ and $x_{\ell+1}$ are identical, and to replace \eqref{eq:mlp-multicouches} with $L$ steps:
\eql{
	x_{\ell+1} = x_\ell + \frac{1}{L} U_\ell^\top \sigma( V_\ell x_\ell + b_\ell), \label{eq:resnet}
}
where $U_\ell, V_\ell \in \mathbb{R}^{n \times d}$ are weight matrices, with $n$ being the number of neurons per layer (which, as in \eqref{eq:mlp-2-couches}, can be increased to enlarge the function space).
The intuition for the success of~\eqref{eq:resnet} is that this formula allows, unlike~\eqref{eq:mlp-multicouches}, for steps that are small deformations near identity mappings. This makes the network $f_\theta(x_0) = x_L$ obtained after $L$ steps well-posed even when $L$ is large, and it can be rigorously proven~\cite{marion2023implicit} that this well-posedness is preserved during optimization via gradient descent~\eqref{eq:grad-desc}.

\paragraph{Neural Differential Equation.}
As $L$ approaches $+\infty$, \eqref{eq:resnet} can be interpreted as a discretization of an ordinary differential equation:
\eql{
	\frac{\mathrm{d} x_s}{\mathrm{d} s} = U_s^\top \sigma(V_s x_s + b_s), \label{eq:neuralode}
}
where $s \in [0,1]$ indexes the network depth.
The network $f_\theta(x_0) := x_1$ maps the initialization $x_{s=0}$ to the solution $x_{s=1}$ of \eqref{eq:neuralode} at time $s=1$. The parameters are $\theta= (U_s, V_s, b_s)_{s \in [0,1]}$.
This formalization, referred to as a \textit{neural ODE}, was initially introduced in \cite{chen2018neural} to leverage tools from adjoint equation theory to reduce the memory cost of computing the gradient $\nabla E(\theta)$ during backpropagation. It also establishes a connection to control theory, as training via gradient descent \eqref{eq:grad-desc} computes an \textit{optimal control} $\theta$ that interpolates between the data $x^i$ and labels $y^i$. 
However, the specificity of learning theory compared to control theory lies in the goal of computing such control using gradient descent \eqref{eq:grad-desc}. To date, no detailed results exist on this. Nonetheless, in his thesis, Raphael Barboni \cite{barboni2024understanding} demonstrated that if the network is initialized near an interpolating network, gradient descent converges to it.

\section{Generative AI for Vector Data} \label{sec:ia-gen}

\paragraph{Self-Supervised Pre-Training.}

The remarkable success of large generative AI models for vector data, such as images (often referred to as ``diffusion models''), and for text (large language models or LLMs), relies on the use of new large-scale network architectures and the creation of new training tasks known as \guill{self-supervised}.
Manually defining the labels $y^i$ through human intervention is too costly, so these are calculated automatically by solving simple tasks. For images, these involve denoising tasks, while for text, they involve predicting the next word.
A key advantage of these simple tasks (known as pre-training tasks) is that it becomes possible to use a pre-trained network $f_\theta$ generatively. Starting from an image composed of pure noise and iterating on the network, one can randomly generate a realistic image~\cite{sohl2015deep}. Similarly, for text, starting from a prompt and sequentially predicting the next word, it is possible to generate text, for example, to answer questions~\cite{brown2020language}.
We will first describe the case of vector data generation, such as images, and address LLMs for text in the following section.

\paragraph{Sampling as an Optimization Problem.}

Generative AI models aim to generate (or \guill{sample}) random vectors $x$ according to a distribution $\beta$, which is learned from a large training dataset $(x^i)_i$. This distribution is obtained by \guill{pushing forward} a reference distribution $\alpha$ (most commonly an isotropic Gaussian distribution $\alpha = \mathcal{N}(0, \text{Id})$) through a neural network $f_\theta$. Specifically, if $X \sim \alpha$ is distributed according to $\alpha$, then $f_\theta(X) \sim \beta$ follows the law $\beta$, which is denoted as $(f_\theta)_\sharp \alpha = \beta$, where $_\sharp$ represents the pushforward operator (also known as the image measure in probability theory).

Early approaches to generative AI, such as \guill{Generative Adversarial Networks} (GANs)~\cite{goodfellow2020generative}, attempted to directly optimize a distance $D$ between probability distributions (e.g., an optimal transport distance~\cite{peyre2019computational}):
\eql{\label{eq:fit-sampling}
	\min_\theta D( (f_\theta)_\sharp \alpha, \beta ).
}
This problem is challenging to optimize because computing $D$ is expensive, and $\beta$ must be approximated from the data $(x^i)_i$.

\paragraph{Flow-Based Generation.}

Recent successes, particularly in the generation of vector data, involve neural networks $f_\theta$ that are computed iteratively by integrating a flow~\cite{papamakarios2021normalizing}, similar to a neural differential equation~\eqref{eq:neuralode}:
\eql{\label{eq:neural-sampling}
	f_\theta(x_0) := x_1
	\quad\text{where}\quad
	\frac{\mathrm{d} x_s}{\mathrm{d} s} = g_\theta(x_s, s),
}
where $g_\theta : (x, s) \in \mathbb{R}^d \times \mathbb{R} \to \mathbb{R}^d$. The input space includes an additional temporal dimension $s$, compared to $f_\theta$. The most effective neural networks for images are notably U-Nets~\cite{ronneberger2015u}.

The central mathematical question is how to replace~\eqref{eq:fit-sampling} with a simpler optimization problem when $f_\theta$ is computed by integrating a flow~\eqref{eq:neural-sampling}. An extremely elegant solution was first proposed in the context of diffusion models~\cite{sohl2015deep} (corresponding to the specific case $\alpha = \mathcal{N}(0, \text{Id})$) and later generalized under the name \guill{flow matching}~\cite{lipman2022flow}.

The main idea is that if $x_0 \sim \alpha$ is initialized randomly, then $x_s$, the solution at time $s$ of \eqref{eq:neural-sampling}, follows a distribution $\alpha_s$ that interpolates between $\alpha_0 = \alpha$ and $\alpha_1$, which is desired to match $\beta$. This distribution satisfies a conservation equation:
\eql{\label{eq:conservation-law}
	\partial_s \alpha_s + \mathrm{div}( \alpha_s v_s ) = 0,
}
where the vector field is defined as $v_s(x) := g_\theta(x, s)$.
This equation is similar to the evolution of neuron distributions during optimization~\eqref{eq:edp-mlps}, but it is simpler because $v_s$ does not depend on the distribution $\alpha_s$, making the equation linear in $\alpha_s$.

\paragraph{Denoising Pre-Training.}

The question is how to find a vector field $v_s(x) = g_\theta(x, s)$ such that $\alpha_1 = \beta$, i.e., the final distribution matches the desired one.
If $v_s$ is known, the distribution $\alpha_s$ is uniquely determined. The key idea is to reason in reverse: starting from an interpolation $\alpha_s$ satisfying $\alpha_0 = \alpha$ and $\alpha_1 = \beta$, how can we compute a vector field $v_s$ such that~\eqref{eq:conservation-law} holds? There are infinitely many possible solutions (since $v_s$ can be modified by a conservative field), but for certain specific interpolations $\alpha_s$, remarkably simple expressions can be derived.

One example is the interpolation obtained via barycentric averaging: take a pair $x_0 \sim \alpha$ and $x_1 \sim \beta$, and define $\alpha_s$ as the distribution of $(1-s)x_0 + sx_1$. It can be shown~\cite{lipman2022flow} that an admissible vector field is given by a simple conditional expectation:
\eql{\label{eq:denoising-diffusion}
	v_s(x) = \mathbb{E}_{x_0 \sim \alpha, x_1 \sim \beta}\big( x_1 - x_0 \mid (1-s)x_0 + sx_1 = x \big).
}
The key advantage of this formula is that the conditional mean corresponds to a linear regression, which can be approximated using a neural network $v_s \approx g_\theta(\cdot, s)$. The expectation over $x_1 \sim \beta$ can be replaced by a sum over training data $(x^i)_i$, leading to the following optimization problem:
$$
	\min_{\theta} E(\theta) := \int_0^1 \mathbb{E}_{x_0 \sim \alpha} \sum_i \| x^i - x_0 - g_\theta((1-s)x_0 + sx^i, s) \|^2 \, \mathrm{d}s.
$$
This function $E(\theta)$ is an empirical risk function~\eqref{eq:erm}, and similar optimization techniques are used to find an optimal $\theta$. In particular, stochastic gradient descent efficiently handles both the integral $\int_0^1$ and the expectation $\mathbb{E}_{x_0 \sim \alpha}$.

The problem~\eqref{eq:denoising-diffusion} corresponds to an unsupervised pre-training task: there are no labels $y^i$, but an artificial supervision task is created by adding random noise $x_0$ to the data $x^i$. This task is called denoising.
We will now see that, for textual data, a different pre-training task is used.

\section{Generative AI for Text} \label{sec:ia-gen-texte}

\paragraph{Tokenization and Next-Word Prediction.}

Generative AI methods for text differ from those used for vector data generation. The neural network architectures are different (they involve transformers, as we will describe), and the pre-training method is based not on denoising but on next-word prediction~\cite{sutskever2014sequence}.
It is worth noting that these transformer neural networks are now also used for image generation~\cite{dosovitskiy2020image}, but specific aspects related to the causal nature of text remain crucial.
The first preliminary step, called \guill{tokenization}, consists of transforming the input text into a sequence of vectors $X = (x[1], \ldots, x[P])$, where the number $P$ is variable (and may increase, for instance, when generating text from a prompt). Each token $x[p]$ is a vector that encodes a group of letters, generally at the level of a syllable.
The neural network $x = f_\theta(X)$ is applied to all the tokens and predicts a new token $x$. During training, a large corpus of text $(X^i)_i$ is available. If we denote by $\tilde X^i$ the text $X^i$ with the last token $x^i$ removed, we minimize an empirical risk for prediction:
$$
	\min_{\theta} E(\theta) := \sum_i \ell( f_\theta(\tilde X^i), x^i ),
$$
which is exactly similar to~\eqref{eq:erm}.
When using a pre-trained network $f_\theta$ for text generation, one starts with a prompt $X$ and iteratively adds a new token in an auto-regressive manner: $X \leftarrow [X, f_\theta(X)]$.

\paragraph{Transformers and Attention.}

The large networks used for text generation tasks are Transformer networks~\cite{vaswani2017attention}. Unlike ResNet~\eqref{eq:resnet}, these networks $f_\theta$ no longer operate on a single vector $x$, but on a set of vectors $X = (x[1], \ldots, x[P])$ of size $P$.
In Transformers, the ResNet layer~\eqref{eq:resnet}, operating on a single vector, is replaced by an \textit{attention} layer, where all tokens interact through a barycentric combination of vectors $(V x[q])_q$ where $V$ is a matrix:
\eql{
    A_{\omega}(X)_p := \sum_q M_{p,q} (Vx[q]), \quad \text{where} \quad M_{p,q} := \frac{e^{\langle Q  x[p], K  x[q] \rangle}}{\sum_\ell e^{\langle Q  x[p], K  x[\ell] \rangle}} \quad\text{with}\quad \omega := (Q,K,V). \label{eq:attention}
}
The coefficients $M_{p,q}$ are computed by normalized correlation between $x[p]$ and $x[q]$, depending on two parameter matrices $Q$ and $K$. A Transformer $f_\theta(X)$ is then defined similarly to ResNet~\eqref{eq:resnet}, iterating $L$ attention layers with residual connections:
\eql{
    X_{\ell+1} = X_\ell + \frac{1}{L} A_{\omega_\ell}(X_\ell). \label{eq:transformers}
}
The parameters $\theta$ of the Transformer $f_\theta(X_0)=X_L$ with $L$ layers are $\theta = (\omega_\ell = (Q_\ell, K_\ell, V_\ell))_{0}^{L-1}$. 

This description is simplified: in practice, a Transformer network also integrates normalization layers and MLPs operating independently on each token. For text applications, attention must be causal, imposing $M_{i,j} = 0$ for $j > i$. This constraint is essential for recursively generating text and ensuring the next-word prediction task is meaningful.

\paragraph{Mean-Field Representation of Attention.}

The attention mechanism \eqref{eq:attention} can be viewed as a system of interacting particles, and \eqref{eq:transformers} describes the evolution of tokens across depth. As with ResNet~\eqref{eq:neuralode}, in the limit $L \to +\infty$, one can consider a system of coupled ordinary differential equations:
\eql{
  \frac{\mathrm{d} X_s}{\mathrm{d} s}  =  A_{\omega_s}(X_s). \label{eq:transformer-ode}
}
A crucial point is that for non-causal Transformers, the system $X_s = (x_s[p])_p$ is invariant under permutations of the indices. Thus, this system can be represented as a probability distribution $\mu_s := \frac{1}{P} \sum_p \delta_{x_s[p]}$ over the token space. This perspective was adopted in Michael Sander's thesis \cite{sander2022sinkformers}, which rewrites attention as:
\[
	\mathcal{A}_{\omega}(x) := \frac{\int e^{\langle Q x, K  x' \rangle} V x' \, \mathrm{d}\mu(x')}{\int e^{\langle Q x, K x' \rangle} \, \mathrm{d}\mu(x')}
	 \quad\text{where}\quad \omega := (Q,K,V).
\]
The particle system \eqref{eq:transformer-ode} then becomes a conservation equation for advection by the vector field $\mathcal{A}_{\omega_s}(\mu)$:
\eql{
    \partial_s \mu_s + \mathrm{div}(\mu_s  \mathcal{A}_{\omega_s}(\mu_s)) = 0. \label{eq:edp-transformer}
}
Surprisingly, this yields a McKean-Vlasov-type equation, similar to the one describing the training of two-layer MLPs \eqref{eq:edp-mlps}, but with the velocity field $\mathcal{A}_{\omega}(\mu)(x)$ replacing $\mathcal{V}(\rho)(u,v,b)$. 
However, here the evolution occurs in the token space $x$ rather than in the neuron space $(u,v,b)$, and the evolution variable corresponds to depth $s$, not the optimization time $t$. 

Unlike \eqref{eq:edp-mlps}, the evolution \eqref{eq:edp-transformer} is not a gradient flow in the Wasserstein metric \cite{sander2022sinkformers}. Nonetheless, in certain cases, this evolution can be analyzed, and it can be shown that the measure $\mu_s$ converges to a single Dirac mass \cite{geshkovski2024emergence} as $s \to +\infty$: the tokens tend to cluster. Better understanding these evolutions, as well as the optimization of parameters $\theta = (Q_s, K_s, V_s)_{s \in [0,1]}$ via gradient descent \eqref{eq:grad-desc}, remains an open problem. This problem can be viewed as a control problem for the PDE \eqref{eq:edp-transformer}.

\section*{Conclusion}

Mathematics plays a critical role in understanding and improving the performance of deep architectures while presenting new theoretical challenges. The emergence of Transformers and generative AI raises immense mathematical problems, particularly for better understanding the training of these networks and exploring the structure of "optimal" networks. One essential question remains: whether an LLM merely interpolates training data or is capable of genuine reasoning. Moreover, issues of resource efficiency and privacy in AI system development demand significant theoretical advancements, where mathematics will play a pivotal role. Whether designing resource-efficient models, ensuring compliance with ethical standards, or exploring the fundamental limits of these systems, mathematics is poised to be an indispensable tool for the future of artificial intelligence.

\bibliographystyle{plain}
\bibliography{bibliography}

\end{document}